\documentclass[12pt, reqno]{amsart}
\usepackage{amsmath, amsthm, amscd, amsfonts, amssymb, graphicx, color}
\usepackage[bookmarksnumbered, colorlinks, plainpages]{hyperref}
\usepackage{subfig}

\usepackage{tikz}
\usetikzlibrary{arrows.meta}
 \usetikzlibrary{calc}
\usepackage[T1,T2A]{fontenc}
\newtheorem{theorem}{Theorem}[section]
\newtheorem{lemma}[theorem]{Lemma}
\newtheorem{proposition}[theorem]{Proposition}
\newtheorem{corollary}[theorem]{Corollary}
\theoremstyle{definition}
\newtheorem{definition}[theorem]{Definition}
\newtheorem{example}[theorem]{Example}

\newtheorem{conjecture}[theorem]{Conjecture}

\newtheorem{problem}[theorem]{Problem}
\theoremstyle{remark}
\newtheorem{remark}[theorem]{Remark}
\numberwithin{equation}{section}

\begin{document}
\setcounter{page}{1}

\centerline{}

\centerline{}

\title[Permanent Saturation of Metric Graphs]{Upper Bound for Permanent Saturation of Metric Graphs using Interval Exchange Transformations}

\author{Ermolaev Egor$^{1}$}
\thanks{$^{1}$ Faculty of Computer Science, HSE University, Moscow, Russia}
\author{Chernyshev Vsevolod$^{2}$}
\thanks{$^{2}$ Ulm University, Ulm, Germany}

\author{Skripchenko Alexandra$^{3}$}
\thanks{$^{3}$ International Laboratory of Cluster Geometry, HSE University, Moscow, Russia}

\begin{abstract} We study upper bounds for the moment of permanent $\varepsilon$-saturation in finite metric graphs. The dynamics is generated by moving points travelling with unit speed along edges and branching into all outgoing directions whenever they reach a vertex. We first reformulate this branched dynamics in terms of birth times at vertices and prove a sufficient same-time criterion for permanent $\varepsilon$-saturation. The main rigorous estimate is obtained from a rotation, regarded as a two-interval exchange transformation. More precisely, if the graph contains two closed walks based at the initial vertex whose lengths have irrational ratio, then the covering properties of the corresponding rotation imply an explicit upper bound for the permanent saturation time. In particular, bounded-type rotations yield a bound of order $\varepsilon^{-1}$. We also construct a more general auxiliary interval exchange transformation on the set of oriented edges. This construction depends on cyclic orders at the vertices and organizes the ordered edge-state data of the graph. Since the branched graph dynamics is non-invertible, whereas an interval exchange transformation is invertible away from discontinuities, this auxiliary IET is not identified with the full graph dynamics. Instead, we formulate the additional birth-time transfer property required for recurrence estimates of the auxiliary IET to imply saturation bounds. We also discuss rotation-type and non-rotation examples of graph-induced self-similar IETs, together with numerical illustrations for star and complete graphs. \end{abstract}

\maketitle
\tableofcontents

\section{Introduction}

Metric graphs provide a natural framework for studying dynamical processes on
networks. A metric graph is a graph whose edges are assigned positive lengths,
so that each edge is identified with a real interval. Such objects appear in
several areas of mathematics, including tropical geometry and combinatorics
\cite{KralHl}, as well as in transport-type problems on networks
\cite{Cher4,Cher5,Cher6,Cher7}. In the present paper we study a deterministic
branched dynamics on a finite metric graph: moving points travel along edges
with unit speed, and whenever a point reaches a vertex, new copies are generated
in all outgoing directions. Equivalently, one may regard the dynamics as the
motion of centers of expanding intervals, but throughout the paper we use the
language of moving points.

The main quantity of interest is the moment of permanent
$\varepsilon$-saturation. Informally, this is the first time after which the
set of moving points remains an $\varepsilon$-net of the whole metric graph
forever. This problem was introduced by Eliseev and Chernyshev in
\cite{CherA}, where upper bounds were obtained under incommensurability
assumptions on the edge lengths. The present work continues this line of
research and develops an interval-exchange viewpoint on the problem.

Interval exchange transformations are bijections of an interval obtained by
cutting it into finitely many subintervals and rearranging them by translations.
They form a classical family of low-complexity dynamical systems and have been
studied extensively from the viewpoints of minimality, ergodicity,
Rauzy--Veech induction, and Teichmüller dynamics; see, for example,
\cite{Keane,Ma,Ve1,Ve2,Vian}. In this paper, interval exchange transformations
enter in two related but distinct ways.

First, rotations are used as two-interval exchange transformations. This gives
the main rigorous saturation estimate of the paper. The key observation is that
birth times at vertices can be described in terms of lengths of walks. If the
graph contains two closed walks based at the initial vertex with lengths $A$
and $B$ such that $A/B\notin\mathbb Q$, then the semigroup
$\{nA+mB:n,m\ge0\}$ produces sufficiently dense birth times. This density is
controlled by the covering time of the rotation
$$
R_{A,B}(x)=x+A\pmod B.
$$
We prove that the corresponding rotation covering estimate implies an explicit
upper bound for the moment of permanent $\varepsilon$-saturation. In the
bounded-type case, this gives a bound of order $\varepsilon^{-1}$.

Second, we construct a more general auxiliary interval exchange transformation
associated with ordered graph data. For each vertex one fixes a cyclic order on
the outgoing oriented edges, and one also chooses a total order on the set of
oriented edges. These data define an auxiliary IET whose intervals are labelled
by oriented edges and whose lengths are inherited from the metric graph. This
construction is useful for organizing the edge-state combinatorics and for
formulating possible Rauzy--Veech-type refinements.

It is important, however, that this auxiliary IET should not be identified
directly with the full branched graph dynamics. The graph dynamics is
non-invertible because of branching, while an IET is invertible away from its
discontinuities. Therefore recurrence estimates for the auxiliary IET imply
saturation bounds only after an additional transfer statement has been proved,
namely a statement converting IET recurrence into density of birth times at the
vertices. We formulate this IET-to-birth-time transfer property explicitly and
derive a conditional Rauzy--Veech saturation estimate under this assumption.

This distinction is one of the main structural points of the paper. The
rotation argument gives an unconditional saturation bound, while the general
auxiliary-IET construction provides a framework for further refinements. In
particular, non-rotation auxiliary IETs can arise from suitable choices of
cyclic orders and top orders. We discuss small self-similar examples, including
rotation-type examples with non-zero Sah--Arnoux--Fathi invariant and a small
non-rotation example. Finding larger non-rotation graph-induced self-similar
examples with non-zero SAF invariant appears to be an interesting separate
problem.

Thus the paper has two complementary aims. On the one hand, it gives a rigorous
birth-time formulation of the branched graph dynamics and proves an explicit
rotation-based upper bound for permanent saturation. On the other hand, it
constructs an auxiliary oriented-edge IET associated with ordered graph data and
clarifies the additional transfer mechanism that would be needed in order to
turn general Rauzy--Veech recurrence estimates into saturation estimates for
the original graph dynamics. This separates the unconditional part of the
argument from the conditional interval-exchange framework and avoids identifying
the invertible auxiliary IET with the non-invertible branched dynamics.

The paper is organized as follows. Section~2 recalls the necessary background
on metric graphs, interval exchange transformations, rotations, and
Rauzy--Veech induction. Section~3 reformulates the graph dynamics in terms of
birth times, proves the main rotation-based saturation estimate, constructs the
auxiliary oriented-edge IET, and states the conditional Rauzy--Veech refinement.
Section~4 discusses self-similar graph-induced IETs and the role of the
Sah--Arnoux--Fathi invariant. Section~5 contains numerical illustrations for
star graphs and complete graphs, together with a brief discussion of possible
applications. Section~6 discusses limitations, conjectural extensions, and open
problems.

\section{Preliminaries}

In this section we recall the notions from the theory of interval exchange transformations and metric graph dynamics that will be used in the sequel. We first review interval exchange transformations and Rauzy--Veech induction, then recall recurrence estimates relevant to our argument, and finally summarize the definitions concerning metric graphs, moving points, and saturation introduced in \cite{CherA}.

\subsection{Interval exchange transformations}

Let $I=[0,L)$ be an interval, and let
\[
0=a_0<a_1<\dots<a_d=L
\]
be a partition of $I$ into subintervals
\[
I_\alpha=[a_{\alpha-1},a_\alpha), \qquad \alpha=1,\dots,d.
\]
Let $\lambda_\alpha=|I_\alpha|$ and write
\[
\lambda=(\lambda_1,\dots,\lambda_d).
\]
An \emph{interval exchange transformation} (IET) is a bijection
\[
T:I\to I
\]
such that, for each $\alpha=1,\dots,d$, the restriction of $T$ to $I_\alpha$ is a translation, and the images $T(I_\alpha)$ form a partition of $I$ into subintervals arranged according to a permutation.

Equivalently, an IET is determined by combinatorial data $\pi=(\pi_0,\pi_1)$ and a length vector $\lambda$, where $\pi_0$ and $\pi_1$ describe the order of the subintervals before and after the exchange. We write $T=T_{\pi,\lambda}$ for the corresponding transformation.

\subsection{Rauzy--Veech induction and Rauzy classes}

We briefly recall the classical Rauzy--Veech induction for interval exchange transformations; see \cite{Ra,Ve2,Vian,AvGoYo}. Let $T=T_{\pi,\lambda}$ be an IET on $d$ subintervals. Comparing the lengths of the rightmost exchanged intervals in the top and bottom orders determines the type of $T$, denoted by $\varepsilon\in\{0,1\}$. Removing the shorter of these two intervals and taking the first return map to the remaining interval produces a new IET, denoted by $\mathcal R(T)$.

This procedure defines the Rauzy--Veech induction, which updates both the combinatorial data and the length vector. Iterating it yields a sequence of IETs together with Rauzy--Veech matrices that encode the corresponding renormalization dynamics.

The connected components of the directed graph generated by Rauzy moves on irreducible permutations are called \emph{Rauzy classes}; see, for example, \cite{Ra,Ve2,Vian}.
\subsection{Suspensions and translation surfaces}\label{subsec:susp}

We briefly recall the suspension of an interval exchange transformation, which
is the standard topological model used in the theory of IETs; see, for example,
\cite{Ve2,ViaT,AvGoYo,Vian}.

\begin{figure}[t]
  \centering
  \includegraphics[width=.6\textwidth]{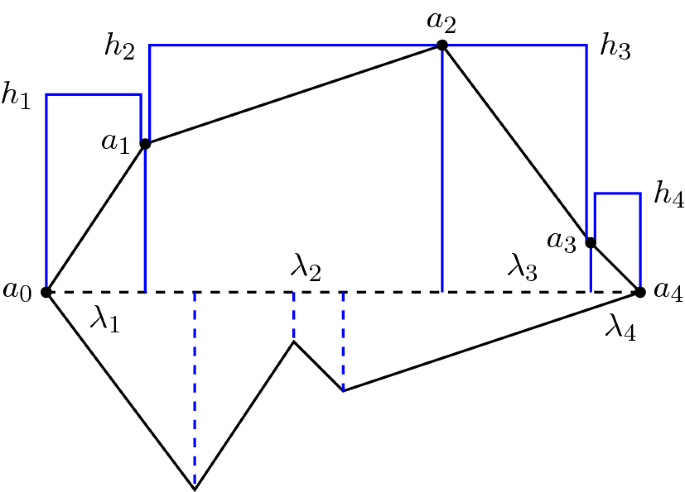}
  \caption{%
  Polygonal construction of the suspension over an IET.
  The base interval $[a_0,a_4)$ is partitioned at $a_1,a_2,a_3$ into pieces of
  lengths $\lambda_1,\lambda_2,\lambda_3,\lambda_4$ (bottom dashed line).
  Vertical segments of heights $h_1,\dots,h_4$ specify the roof; joining their
  tops gives the upper broken boundary. Gluing, by translations, the horizontal
  copies of each subinterval on the top and bottom in the orders $\pi_0$ and
  $\pi_1$ yields a translation surface, and the vertical straight-line flow
  there is the special flow over the IET with roof $h$.}
  \label{fig:polygon_suspension}
\end{figure}

Let $T=T_{\pi,\lambda}\colon[0,L)\to[0,L)$ be an IET with irreducible
permutation data $(\pi_0,\pi_1)$, and fix a roof function $r\colon[0,L)\to
\mathbb R_{>0}$ which is constant on each top subinterval, $r(x)\equiv h_i$ on
$I_i$. The \emph{suspension space} is
\[
X_{T,r}=\bigl\{(x,s):x\in[0,L),\ 0\le s<r(x)\bigr\}\big/\!\sim,\qquad
(x,r(x))\sim (Tx,0),
\]
and the \emph{suspension flow} is the unit-speed vertical flow on $X_{T,r}$
with identifications at the roof. Equivalently, one may form rectangles
$R_i=[0,\lambda_i]\times[0,h_i]$ and glue their vertical sides by translations
according to $(\pi_0,\pi_1)$; this is the zippered-rectangles construction of a
translation surface $S(\pi,\lambda,h)$. The polygonal version of the same
construction, illustrated below, builds the same surface from a single polygon
whose lower side is the base interval and whose upper broken side is determined
by the roof; see Figure~\ref{fig:polygon_suspension}. In each picture, the
vertical straight-line flow is measurably isomorphic to the suspension flow
over $T$, and $T$ is the first-return map to the base interval.

The Teichmüller flow $g_t=\mathrm{diag}(e^t,e^{-t})$ acts on translation
charts by stretching the horizontal direction and contracting the vertical
one; Rauzy--Veech induction is a discrete cross-section of $(g_t)$
\cite{ViaT,AvGoYo,Ve2}. In the present paper this picture is used only as
background for the conjectural Lyapunov-type refinements of
Section~\ref{sec:discussion}. The rigorous saturation estimate of
Theorem~\ref{main} is proved through birth times at vertices, not through an
identification of the full branched graph dynamics with a suspension flow:
indeed, the branched dynamics is non-invertible because of branching, whereas
an IET and its suspension are invertible away from singularities.

\subsection{Recurrence and return times.} 
Let $\mu$ be a probability measure on a metric space $(X, d)$ defined on Borel subsets of $X$ and let $T$ be a measure preserving transformation. Given a measurable subset $E$ with positive measure and a point $x \in E$ one can consider the first return time $R_E(x)=\inf \left\{j \geqslant 1 \mid T^j(x) \in E\right\}$ of the point $x$ in the set $E$. By the Poincaré recurrence theorem \cite{Vian} one has $R_E(x)<+\infty$ for almost every $x$ and Kac's lemma \cite{Kac} states that $\int_E R_E(x) \mathrm{d} \mu \leqslant 1$ where equality holds for ergodic transformations.
Let $r>0$ and let $\tau_r(x)$ be the return time to the $r$-neighbourhood of $x$
$$
\tau_r(x)=\min \left\{j \geqslant 1: d\left(T^j x, x\right)<r\right\}
$$
and let $\tau_r(x, y)$ be the hitting time of an $r$-neighbourhood
$$
\tau_r(x, y)=\min \left\{j \geqslant 1: d\left(T^j x, y\right)<r\right\} .
$$

For almost every interval exchange transformation we have that \cite{Bosh, KimS}
$$
\begin{array}{ll}
\lim _{r \rightarrow 0} \frac{\log \tau_r(x)}{-\log r}=1 & \text { a.e. } x, \\
\lim _{r \rightarrow 0} \frac{\log \tau_r(x, y)}{-\log r}=1 & \text { a.e. } x .
\end{array}
$$

\subsection{Metric graphs and saturation dynamics}

In this subsection we recall the definitions from \cite{CherA} that will be used throughout the paper.
\begin{definition}
A metric graph is a triple $G=(V,E,l)$, where $(V,E)$ is a graph and
\[
l:E\to\mathbb R_{>0}
\]
is a length function. Each edge $e\in E$ is identified with an interval of length $l(e)$.
\end{definition}

Throughout the paper we restrict attention to undirected simple metric graphs. Thus there is at most one edge between any two vertices and no loops. When convenient, an edge joining vertices $v_i$ and $v_j$ may be denoted by $e_{ij}$.

\begin{definition} 
Let $G=(V,E,l)$ be a metric graph, $v^* \in V$. Dynamical system of dispersing moving points $\mathcal{W}\left(G, v^*\right)$ is defined by the following rules:
\begin{enumerate}
    \item Time and space are continuous, with $\mathbb{R}_{\geq 0}$ and $V \cup\left(\bigsqcup_{e \in E} e\right)$ being their respective domains.
    \item Each moving point lies on an oriented edge and moves at unit speed toward the head vertex of that oriented edge. Its direction changes only when it reaches a vertex, according to the branching rule below.
    \item When a moving point reaches a vertex $v$, the point is removed and new
moving points are created on all oriented edges issuing from $v$. In particular,
if the point arrived at $v$ along an oriented edge $a$, then one of the new
copies is created on the reverse oriented edge $\bar a$. Each new point starts
at $v$ and moves away from $v$ at unit speed. The time $t$ at which this occurs
is called a birth time.
    \item At time \(0\), moving points are created at \(v^*\) on all oriented edges issuing from \(v^*\).
\end{enumerate}
\end{definition}

\begin{definition}
Let $G=(V,E,l)$ be a metric graph, $v^* \in V$, $\varepsilon \in \mathbb R_{>0}$, $t\in\mathbb R_{\ge0}$. The system $\mathcal{W}\left(G, v^*\right)$ is said to be $\varepsilon$-saturated at time $t$ if its moving points form an $\varepsilon$-net on $V \cup\left(\bigsqcup_{e \in E} e\right)$ at time $t$.
\end{definition}

\begin{definition}
    Let $G=(V,E,l)$ be a metric graph, $v^* \in V$, $\varepsilon \in \mathbb R_{>0}$. A connected non-empty set $T \subseteq \mathbb{R}_{\geq 0}$ is called an $\varepsilon$-saturation interval of $\mathcal{W}\left(G, v^*\right)$ if $\mathcal{W}\left(G, v^*\right)$ is $\varepsilon$-saturated at any $t \in T$ and if for any $T^{\prime} \subseteq \mathbb{R}_{\geq 0}$ with $T \subset T^{\prime}$ there exists $t^{\prime} \in T^{\prime}$ such that $\mathcal{W}\left(G, v^*\right)$ is not $\varepsilon$-saturated at time $t^{\prime}$. We call $\inf T$ and $\sup T$ an $\varepsilon$-saturation moment and an $\varepsilon$-desaturation moment of $\mathcal{W}\left(G, v^*\right)$ respectively.
\end{definition}

\begin{definition}
     Let $G=(V,E,l)$ be a metric graph, $v^* \in V$, $\varepsilon \in \mathbb R_{>0}$, and let $T$ be an $\varepsilon$-saturation interval of $\mathcal W(G,v^*)$. We call $T$ the interval of permanent $\varepsilon$-saturation, if $\sup T=+\infty$. Likewise, we call $\inf T$ the moment of permanent $\varepsilon$-saturation.
\end{definition}

In what follows, our main quantity of interest is the moment of permanent $\varepsilon$-saturation. For brevity, when $\varepsilon$ is fixed, we sometimes write simply \emph{permanent saturation moment}.

\section{Birth times and interval-exchange models}\label{sec:graph_to_iet}

In this section we relate the branched dynamics on a metric graph to recurrence
problems for interval exchange transformations. We first reformulate the graph
dynamics in terms of birth times at vertices. This gives an sufficient same-time
criterion for permanent saturation. We then explain how recurrence of a rotation,
viewed as a two-interval exchange transformation, yields an explicit saturation
bound. Finally, we describe a more general auxiliary IET associated with ordered
graph data and state precisely what additional transfer property is required in
order to use its recurrence estimates for the full branched dynamics.

\subsection{Birth times and a saturation criterion}
\label{subsec:birth_times}

Let $G=(V,E,l)$ be a finite connected undirected metric graph, and let
$v^*\in V$ be the initial vertex. We denote by $\overrightarrow E$ the set of
oriented edges of $G$. If $a=(u\to w)\in\overrightarrow E$, we write
$$
t(a)=u,\qquad h(a)=w,\qquad l(a)=l(\{u,w\}).
$$

A finite walk from $v^*$ to a vertex $v\in V$ is a finite sequence
$$
\omega=(v_0,v_1,\ldots,v_n)
$$
such that
$$
v_0=v^*,\qquad v_n=v,\qquad \{v_{i-1},v_i\}\in E
\quad\text{for }i=1,\ldots,n.
$$
The metric length of $\omega$ is
$$
|\omega|_l
=
\sum_{i=1}^{n} l(\{v_{i-1},v_i\}).
$$
We also include the empty walk at $v^*$, whose length is $0$.

For $v\in V$, define
$$
\mathcal A_v
=
\{\,|\omega|_l:\omega \text{ is a finite walk from } v^* \text{ to } v\,\}.
$$
Thus $0\in\mathcal A_{v^*}$. The set $\mathcal A_v$ is the set of times at
which moving points are born at the vertex $v$.

\begin{lemma}
\label{lem:arrival_times}
For every vertex $v\in V$, the set $\mathcal A_v$ coincides with the set of
birth times at $v$ in the branched dynamics $\mathcal W(G,v^*)$.
\end{lemma}

\begin{proof}
Let $\omega=(v_0,\ldots,v_n)$ be a finite walk from $v^*$ to $v$. Since the
dynamics creates a copy along every outgoing oriented edge whenever a moving
point reaches a vertex, there exists a descendant following the walk $\omega$.
This descendant reaches $v$ at time $|\omega|_l$. Hence every element of
$\mathcal A_v$ is a birth time at $v$.

Conversely, suppose that a moving point is born at $v$ at time $t$. Tracing
the ancestry of this point backwards gives a finite walk from $v^*$ to $v$.
Since all moving points travel with unit speed, the time $t$ is equal to the
sum of the lengths of the edges in this walk. Therefore $t\in\mathcal A_v$.
\end{proof}

For an oriented edge $a=(u\to w)$, the birth times of moving points travelling
along $a$ are precisely the elements of $\mathcal A_u$. At time $t$, the
positions of such points on $a$, measured as distance from $u=t(a)$, form
the set
$$
X_a(t)
=
\{\,t-\tau:\tau\in\mathcal A_{t(a)},\ 0\le t-\tau\le l(a)\,\}
\subseteq [0,l(a)].
$$
Indeed, a point born at $t(a)$ at time $\tau$ and sent along $a$ is located
at distance $t-\tau$ from $t(a)$, as long as this quantity belongs to
$[0,l(a)]$.

If $e=\{u,w\}$ is an unoriented edge, let
$$
a=(u\to w),\qquad \bar a=(w\to u).
$$
Identifying $e$ with $[0,l(e)]$, where $0$ corresponds to $u$, the set of
occupied positions on $e$ at time $t$ is
$$
X_e(t)
=
X_a(t)
\cup
\{\,l(e)-s:s\in X_{\bar a}(t)\,\}.
$$

\begin{lemma}[Saturation criterion]
\label{lem:birth_saturation_criterion}
Let $\varepsilon>0$. Suppose that there exists $T\ge0$ such that, for every
oriented edge $a\in\overrightarrow E$ and every $t\ge T$, the set
$$
X_a(t)
=
\{\,t-\tau:\tau\in\mathcal A_{t(a)},\ 0\le t-\tau\le l(a)\,\}
$$
is an $\varepsilon$-net in $[0,l(a)]$. Then the branched dynamics
$\mathcal W(G,v^*)$ is $\varepsilon$-saturated for every $t\ge T$. In
particular,
$$
\tau_s(\varepsilon)\le T.
$$
\end{lemma}

\begin{proof}
Fix $t\ge T$. Let $e=\{u,w\}$ be an edge, and identify $e$ with
$[0,l(e)]$, where $0$ corresponds to $u$. Put $a=(u\to w)$. By
assumption, $X_a(t)$ is an $\varepsilon$-net in $[0,l(e)]$. Hence every
point of $e$ lies within distance $\varepsilon$ of a moving point located on
$e$ at time $t$.

The argument applies to every edge of $G$. Since the vertices are endpoints of
edges, the set of all moving points at time $t$ is an $\varepsilon$-net of
the whole metric graph. Thus $\mathcal W(G,v^*)$ is $\varepsilon$-saturated
for every $t\ge T$, and therefore
$$
\tau_s(\varepsilon)\le T.
$$
\end{proof}

\subsection{Closed walks and a rotation IET}
\label{subsec:closed_walk_rotation}

We now give a concrete interval-exchange mechanism which implies the
birth-time density condition in Lemma~\ref{lem:birth_saturation_criterion}.
The interval exchange used here is a rotation, regarded as a two-interval
exchange transformation.

Let $A,B>0$ and assume that
$$
\frac{A}{B}\notin\mathbb Q.
$$
Let
$$
\alpha=A \bmod B,\qquad 0<\alpha<B.
$$
We consider the rotation
$$
R_{A,B}:\mathbb R/B\mathbb Z\to \mathbb R/B\mathbb Z,
\qquad
R_{A,B}(x)=x+A \pmod B.
$$
Equivalently,
$$
R_{A,B}(x)=x+\alpha \pmod B.
$$
As an interval exchange transformation on $[0,B)$, this is the two-interval
exchange obtained by cutting $[0,B)$ into
$$
I_1=[0,B-\alpha),\qquad I_2=[B-\alpha,B),
$$
and defining
$$
R_{A,B}(x)=
\begin{cases}
x+\alpha, & x\in I_1,\\
x+\alpha-B, & x\in I_2.
\end{cases}
$$
Thus the exchanged intervals have lengths
$$
B-\alpha,\qquad \alpha,
$$
and the corresponding permutation is
$$
\begin{pmatrix}
1&2\\
2&1
\end{pmatrix}.
$$

For $\delta>0$, define $N_{A,B}(\delta)$ as follows. For $N\ge1$, let
$$
S_N(A,B)=\{\,nA \bmod B:0\le n\le N\,\}\subset \mathbb R/B\mathbb Z.
$$
Since $A/B$ is irrational, the points of $S_N(A,B)$ are distinct. We say
that $S_N(A,B)$ has circular gaps at most $\delta$ if, after ordering the
points cyclically on the circle $\mathbb R/B\mathbb Z$, every adjacent gap has
length at most $\delta$. Define
$$
N_{A,B}(\delta)
=
\min
\left\{
N\ge1:
S_N(A,B)\text{ has circular gaps at most }\delta
\right\}.
$$
The number $N_{A,B}(\delta)$ is finite for every $\delta>0$, because the
irrational rotation $R_{A,B}$ is minimal.

Equivalently, $N=N_{A,B}(\delta)$ is the least integer such that, for every
$\rho\in[0,B)$, there exists $0\le n\le N$ with
$$
nA\bmod B
\in
[\rho-\delta,\rho]
\quad\text{on the circle } \mathbb R/B\mathbb Z.
$$
Here $[\rho-\delta,\rho]$ denotes the positively oriented circular interval
of length $\delta$ ending at $\rho$.

\begin{lemma}
\label{lem:rotation_semigroup_density}
Let $A,B>0$, $A/B\notin\mathbb Q$, and let
$$
N=N_{A,B}(\delta).
$$
Then, for every
$$
y\ge NA+B,
$$
there exist integers $n,m\ge0$ such that
$$
0\le y-(nA+mB)\le \delta.
$$
\end{lemma}

\begin{proof}
Write
$$
y=kB+\rho,
\qquad
k\in\mathbb Z_{\ge0},
\qquad
\rho\in[0,B).
$$
By the definition of $N=N_{A,B}(\delta)$, there exists $0\le n\le N$ such
that
$$
r:=nA\bmod B
$$
belongs to the circular interval $[\rho-\delta,\rho]$. Write
$$
nA=qB+r,
\qquad
q\in\mathbb Z_{\ge0},
\qquad
r\in[0,B).
$$
Since $y\ge NA+B$ and $n\le N$, we have
$$
kB=y-\rho\ge NA+B-\rho>NA\ge nA.
$$
Hence
$$
kB>qB+r\ge qB,
$$
and therefore $k\ge q+1$.

There are two cases. First suppose that $r\le \rho$. Then
$$
0\le \rho-r\le\delta.
$$
Set
$$
m=k-q.
$$
Since $k\ge q+1$, we have $m\ge0$. Moreover,
$$
nA+mB=qB+r+(k-q)B=kB+r.
$$
Thus
$$
0\le y-(nA+mB)=\rho-r\le\delta.
$$

Now suppose that $r>\rho$. Since $r$ belongs to the circular interval
$[\rho-\delta,\rho]$, we have
$$
0\le B-r+\rho\le\delta.
$$
Set
$$
m=k-q-1.
$$
Again $m\ge0$. Then
$$
nA+mB=qB+r+(k-q-1)B=(k-1)B+r,
$$
and hence
$$
0\le y-(nA+mB)=B-r+\rho\le\delta.
$$
The proof is complete.
\end{proof}

\subsection{The main saturation estimate}
\label{subsec:main_saturation_estimate}

We now apply the elementary rotation estimate from
Lemma~\ref{lem:rotation_semigroup_density} to the branched dynamics on a metric
graph. This gives the main rigorous saturation bound of the paper.

A closed walk based at $v^*$ is a finite walk
$$
C=(v_0,v_1,\ldots,v_n)
$$
such that
$$
v_0=v_n=v^*.
$$
Its metric length is
$$
l(C)=\sum_{i=1}^{n}l(\{v_{i-1},v_i\}).
$$

\begin{theorem}[Main theorem]
\label{main}
Let $G=(V,E,l)$ be a finite connected undirected metric graph, and let
$v^*\in V$. Suppose that there exist two closed walks $C_1$ and $C_2$
based at $v^*$ whose lengths
$$
A=l(C_1),
\qquad
B=l(C_2)
$$
satisfy
$$
\frac{A}{B}\notin\mathbb Q.
$$
For every vertex $u\in V$, fix a path $P_u$ from $v^*$ to $u$, and write
$$
p_u=l(P_u).
$$
Then, for every
$$
0<\varepsilon<\min_{e\in E}l(e),
$$
the moment of permanent $\varepsilon$-saturation satisfies
$$
\tau_s(\varepsilon)
\le
\max_{a=(u\to w)\in\overrightarrow E}
\left(
p_u+A\,N_{A,B}(\varepsilon)+B+l(a)
\right),
$$
where $N_{A,B}(\varepsilon)$ is the covering time of the rotation
$$
R_{A,B}(x)=x+A\pmod B
$$
defined in Subsection~\ref{subsec:closed_walk_rotation}.
\end{theorem}
\begin{remark}
    The hypothesis of Theorem~\ref{main} is non-vacuous. For example, if the
    initial vertex $v^*$ is incident to two distinct edges $e$ and $e'$, then the
    closed walks
    $$
    C_e=(v^*,w,v^*),
    \qquad
    C_{e'}=(v^*,w',v^*)
    $$
    have lengths $2l(e)$ and $2l(e')$. Hence the hypothesis holds whenever
    $l(e)/l(e')\notin\mathbb Q$. More generally, if $C_1$ and $C_2$ are any two
    closed walks based at $v^*$, then the condition
    $$
    \frac{l(C_1)}{l(C_2)}\notin\mathbb Q
    $$
    fails only on a countable union of rational hyperplanes in the edge-length
    parameter space. Thus, for fixed graph topology and fixed closed walks, the
    condition is generic with respect to the edge lengths.
\end{remark}

\begin{proof}
Let
$$
N=N_{A,B}(\varepsilon)
$$
and define
$$
T_\varepsilon
=
\max_{a=(u\to w)\in\overrightarrow E}
\left(
p_u+AN+B+l(a)
\right).
$$
We prove that the graph is $\varepsilon$-saturated for every
$t\ge T_\varepsilon$.

Since $C_1$ and $C_2$ are closed walks based at $v^*$, the branched
dynamics creates moving points at $v^*$ at all times of the form
$$
nA+mB,
\qquad
n,m\in\mathbb Z_{\ge0}.
$$
Indeed, one descendant may traverse $C_1$ exactly $n$ times and then traverse
$C_2$ exactly $m$ times. After returning to $v^*$, another descendant may
follow the fixed path $P_u$ from $v^*$ to $u$. Hence
$$
p_u+nA+mB\in\mathcal A_u
$$
for every $u\in V$ and all $n,m\ge0$.

Fix $t\ge T_\varepsilon$, an oriented edge $a=(u\to w)$, and a coordinate
$$
s\in[0,l(a)]
$$
on $a$, measured from $u$. We first consider
$$
0\le s\le l(a)-\varepsilon.
$$
Set
$$
y=t-s-p_u.
$$
By the definition of $T_\varepsilon$,
$$
y\ge AN+B.
$$
By Lemma~\ref{lem:rotation_semigroup_density}, there exist integers
$n,m\ge0$ such that
$$
0\le y-(nA+mB)\le\varepsilon.
$$
Define
$$
\tau=p_u+nA+mB.
$$
Then $\tau\in\mathcal A_u$. A moving point born at $u$ at time $\tau$ and
sent along the oriented edge $a$ is, at time $t$, at distance
$$
t-\tau
=
s+\bigl(y-(nA+mB)\bigr)
$$
from $u$. Therefore
$$
s\le t-\tau\le s+\varepsilon\le l(a).
$$
Thus this point is still on the edge $a$, and its position is within distance
$\varepsilon$ of the point with coordinate $s$.

It remains to consider the interval
$$
l(a)-\varepsilon\le s\le l(a).
$$
Applying the previous argument to
$$
s_0=l(a)-\varepsilon,
$$
we obtain a moving point on $a$ whose coordinate belongs to
$$
[l(a)-\varepsilon,l(a)].
$$
This point is within distance $\varepsilon$ of every point with coordinate
$s\in[l(a)-\varepsilon,l(a)]$.

Hence, for every oriented edge $a\in\overrightarrow E$ and every
$t\ge T_\varepsilon$, the set $X_a(t)$ is an $\varepsilon$-net in
$[0,l(a)]$. By Lemma~\ref{lem:birth_saturation_criterion}, the branched
dynamics $\mathcal W(G,v^*)$ is $\varepsilon$-saturated for every
$t\ge T_\varepsilon$. Therefore
$$
\tau_s(\varepsilon)\le T_\varepsilon,
$$
which is the desired estimate.
\end{proof}

\begin{corollary}
\label{cor:main_rotation_asymptotic}
In the setting of Theorem~\ref{main}, suppose that there exist constants
$C>0$, $\gamma\ge1$, and $\varepsilon_0>0$ such that
$$
N_{A,B}(\varepsilon)\le C\varepsilon^{-\gamma}
$$
for all $0<\varepsilon<\varepsilon_0$. Then there exists a constant
$K>0$, depending on $G$, $v^*$, $A$, $B$, $C$, and $\gamma$, such
that
$$
\tau_s(\varepsilon)\le K\varepsilon^{-\gamma}
$$
for all sufficiently small $\varepsilon>0$.

In particular, if the rotation $R_{A,B}$ has bounded type, then
$$
\tau_s(\varepsilon)\le K\varepsilon^{-1}
$$
for all sufficiently small $\varepsilon>0$.
\end{corollary}

\begin{proof}
Substituting
$$
N_{A,B}(\varepsilon)\le C\varepsilon^{-\gamma}
$$
into Theorem~\ref{main} gives
$$
\tau_s(\varepsilon)
\le
AC\varepsilon^{-\gamma}
+
\max_{a=(u\to w)\in\overrightarrow E}
\bigl(p_u+B+l(a)\bigr).
$$
The second term is independent of $\varepsilon$, and for sufficiently small
$\varepsilon$ it can be absorbed into the first term after increasing the
constant.

For bounded-type rotations, the standard covering estimate gives
$$
N_{A,B}(\varepsilon)\le C'\varepsilon^{-1}
$$
for sufficiently small $\varepsilon$. This gives the final claim.
\end{proof}

\begin{example}[Star graph]
\label{ex:star_rotation}
Let $G$ be the star graph with center $o$ and leaves
$$
v_1,\ldots,v_k.
$$
Let
$$
e_i=\{o,v_i\},
\qquad
l(e_i)=\ell_i>0.
$$
Assume that
$$
\frac{\ell_1}{\ell_2}\notin\mathbb Q,
$$
and take the initial vertex to be
$$
v^*=o.
$$
Consider the two closed walks
$$
C_1=(o,v_1,o),
\qquad
C_2=(o,v_2,o).
$$
Their lengths are
$$
A=2\ell_1,
\qquad
B=2\ell_2.
$$
Thus
$$
\frac{A}{B}=\frac{\ell_1}{\ell_2}\notin\mathbb Q.
$$
The associated two-interval exchange transformation is the rotation
$$
R_{A,B}:[0,B)\to[0,B),
\qquad
R_{A,B}(x)=x+A\pmod B.
$$
If
$$
\alpha=A\bmod B,
\qquad 0<\alpha<B,
$$
then $R_{A,B}$ is the IET obtained by cutting $[0,B)$ into intervals of
lengths
$$
B-\alpha,\qquad \alpha,
$$
with permutation
$$
\begin{pmatrix}
1&2\\
2&1
\end{pmatrix}.
$$

For the center $o$, we take $P_o$ to be the empty path, so $p_o=0$. For a
leaf $v_i$, we take $P_{v_i}=(o,v_i)$, so $p_{v_i}=\ell_i$. Applying
Theorem~\ref{main}, we obtain
$$
\tau_s(\varepsilon)
\le
A\,N_{A,B}(\varepsilon)+B+2\max_{1\le i\le k}\ell_i.
$$
If the rotation $R_{A,B}$ has bounded type, then
$$
\tau_s(\varepsilon)\le K\varepsilon^{-1}
$$
for all sufficiently small $\varepsilon>0$.
\end{example}

\subsection{An auxiliary IET on oriented edges}
\label{subsec:auxiliary_oriented_edge_iet}

We now describe a more global auxiliary interval exchange transformation
associated with ordered graph data. This construction is useful for organizing
the combinatorics of oriented edges and for formulating possible
Rauzy--Veech-type refinements. It should not, by itself, be identified with the
full branched graph dynamics: the latter is non-invertible because of branching,
whereas an interval exchange transformation is invertible away from its
discontinuities. Thus recurrence estimates for the auxiliary IET imply
saturation bounds only after one proves an additional transfer statement
relating IET recurrence to density of graph birth times.

Let $G=(V,E,l)$ be a finite connected undirected metric graph. Denote by
$\overrightarrow E$ the set of oriented edges of $G$. For
$a=(u\to v)\in\overrightarrow E$, write
$$
\bar a=(v\to u),\qquad t(a)=u,\qquad h(a)=v,
$$
and set
$$
l(a)=l(\{u,v\}).
$$

For each vertex $v\in V$, fix a cyclic order on the oriented edges issuing
from $v$. This gives a cyclic permutation
$$
\sigma_v:
\{a\in\overrightarrow E:t(a)=v\}
\longrightarrow
\{a\in\overrightarrow E:t(a)=v\}.
$$
Define
$$
\Phi:\overrightarrow E\to\overrightarrow E,
\qquad
\Phi(a)=\sigma_{h(a)}(\bar a).
$$

\begin{lemma}
\label{lem:Phi_permutation}
The map $\Phi$ is a permutation of the finite set $\overrightarrow E$.
\end{lemma}

\begin{proof}
The map $a\mapsto \bar a$ is a bijection of $\overrightarrow E$. For each
vertex $v$, the map $\sigma_v$ is a permutation of the oriented edges with
tail $v$. Hence the composition defining $\Phi$ is a bijection of
$\overrightarrow E$.
\end{proof}

We now construct an interval exchange transformation from these ordered data.
Put
$$
d=|\overrightarrow E|.
$$
Choose a total order on $\overrightarrow E$, written as a bijection
$$
\pi_0:\overrightarrow E\to\{1,\ldots,d\}.
$$
The ordered graph data are the tuple
$$
(G,l,\{\sigma_v\}_{v\in V},\pi_0).
$$
Define a second order
$$
\pi_1:\overrightarrow E\to\{1,\ldots,d\}
$$
by
$$
\pi_1(\Phi(a))=\pi_0(a)
\qquad
\text{for all }a\in\overrightarrow E.
$$
Equivalently,
$$
\pi_1(b)=\pi_0(\Phi^{-1}(b)).
$$

Let
$$
\lambda_a=l(a),
\qquad a\in\overrightarrow E,
$$
and set
$$
L=\sum_{a\in\overrightarrow E}\lambda_a.
$$
We define an interval exchange transformation
$$
T_{\mathrm{aux}}=T_{\pi,\lambda}:[0,L)\to[0,L)
$$
as follows. The top partition consists of intervals
$$
I_a=
\left[
\sum_{\pi_0(b)<\pi_0(a)}\lambda_b,
\sum_{\pi_0(b)\le \pi_0(a)}\lambda_b
\right),
\qquad a\in\overrightarrow E.
$$
The bottom partition consists of intervals
$$
J_a=
\left[
\sum_{\pi_1(b)<\pi_1(a)}\lambda_b,
\sum_{\pi_1(b)\le \pi_1(a)}\lambda_b
\right),
\qquad a\in\overrightarrow E.
$$
The map $T_{\mathrm{aux}}$ translates each $I_a$ onto $J_a$. Thus, for
$x\in I_a$,
$$
T_{\mathrm{aux}}(x)
=
x
-
\sum_{\pi_0(b)<\pi_0(a)}\lambda_b
+
\sum_{\pi_1(b)<\pi_1(a)}\lambda_b.
$$
This is a classical interval exchange transformation with alphabet
$\overrightarrow E$, length vector $\lambda$, and permutation data
$(\pi_0,\pi_1)$.

\begin{remark}
The construction above avoids any length mismatch: the interval $I_a$ and
the interval $J_a$ both have length $\lambda_a=l(a)$, and $T_{\mathrm{aux}}$
translates $I_a$ onto $J_a$. The permutation $\Phi$ determines the order of
the labels in the bottom partition, but the IET itself maps the top interval
with label $a$ to the bottom interval with the same label $a$. Therefore the
symbolic itinerary of a point under the top partition should not automatically
be identified with the deterministic successor orbit
$$
a,\Phi(a),\Phi^2(a),\ldots .
$$
For this reason, the auxiliary IET is an ordering and renormalization device,
not a complete model of the branched graph dynamics.
\end{remark}

\begin{remark}
The full graph dynamics also contains multiplicity information: many descendants
may occupy the same oriented edge state at the same time. In the present paper,
this branching information is handled through the birth-time sets
$\mathcal A_v$ introduced in Subsection~\ref{subsec:birth_times}. We do not
use a multiplicity cocycle in the proof of the main saturation bound. A genuine
skew-product model would require an explicit base transformation, an explicit
fiber map, and a proof that the resulting system factors onto the branched graph
dynamics.
\end{remark}

\begin{example}[Auxiliary IET for a three-armed star graph]
\label{ex:aux_star}
Let $G$ be the star graph with center $o$ and leaves
$v_1,v_2,v_3$. Write
$$
e_i=\{o,v_i\},
\qquad
l(e_i)=\ell_i>0,
\qquad
i=1,2,3.
$$
The oriented-edge alphabet is
$$
\overrightarrow E=\{a_1,a_2,a_3,b_1,b_2,b_3\},
$$
where
$$
a_i=(o\to v_i),
\qquad
b_i=(v_i\to o),
\qquad
i=1,2,3.
$$

Fix the cyclic order at the center by
$$
a_1\mapsto a_2\mapsto a_3\mapsto a_1.
$$
At each leaf $v_i$, the cyclic order is trivial. Hence
$$
\Phi(a_i)=b_i,
\qquad i=1,2,3,
$$
and
$$
\Phi(b_1)=a_2,\qquad
\Phi(b_2)=a_3,\qquad
\Phi(b_3)=a_1.
$$
Thus
$$
\Phi=(a_1\, b_1\, a_2\, b_2\, a_3\, b_3).
$$

Choose the top order
$$
\pi_0=(a_1,b_1,a_2,b_2,a_3,b_3).
$$
Then the bottom order defined by $\pi_1(\Phi(a))=\pi_0(a)$ is
$$
\pi_1=(b_1,a_2,b_2,a_3,b_3,a_1).
$$
The length vector is
$$
\lambda_{a_i}=\lambda_{b_i}=\ell_i,
\qquad i=1,2,3.
$$
Therefore $T_{\mathrm{aux}}$ acts on an interval of total length
$$
L=2(\ell_1+\ell_2+\ell_3),
$$
partitioned into six intervals with lengths
$$
(\ell_1,\ell_1,\ell_2,\ell_2,\ell_3,\ell_3)
$$
in the top order $\pi_0$, and rearranged in the bottom order $\pi_1$.

In this example the bottom order is obtained from the top order by a cyclic
shift. Hence the refined IET is the rotation
$$
T_{\mathrm{aux}}(x)=x-\ell_1 \pmod L.
$$
The refined presentation has removable discontinuities. After merging them,
one obtains the usual rotation of the circle of length $L$. In particular,
the reduced map is minimal and uniquely ergodic whenever
$$
\frac{\ell_1}{L}\notin\mathbb Q.
$$
\end{example}

\begin{example}[Auxiliary IET for a triangle graph]
\label{ex:aux_triangle}
Let $G$ be the triangle with vertices $A,B,C$ and edges
$$
e_{AB}=\{A,B\},
\qquad
e_{BC}=\{B,C\},
\qquad
e_{CA}=\{C,A\}.
$$
Let
$$
a_1=(A\to B),\qquad a_2=(B\to C),\qquad a_3=(C\to A),
$$
and
$$
b_1=(B\to A),\qquad b_2=(C\to B),\qquad b_3=(A\to C).
$$
Choose the cyclic orders
$$
\sigma_A=(a_1\ b_3),\qquad
\sigma_B=(a_2\ b_1),\qquad
\sigma_C=(a_3\ b_2).
$$
Then
$$
\Phi(a_1)=a_2,\qquad
\Phi(a_2)=a_3,\qquad
\Phi(a_3)=a_1,
$$
and
$$
\Phi(b_1)=b_3,\qquad
\Phi(b_3)=b_2,\qquad
\Phi(b_2)=b_1.
$$
Thus
$$
\Phi=(a_1\,a_2\,a_3)(b_1\,b_3\,b_2).
$$

If
$$
\pi_0=(a_1,a_2,a_3,b_1,b_2,b_3),
$$
then
$$
\pi_1=(a_2,a_3,a_1,b_3,b_1,b_2).
$$
The associated length vector is
$$
\lambda=
\bigl(
l(e_{AB}),\,l(e_{BC}),\,l(e_{CA}),\,
l(e_{AB}),\,l(e_{BC}),\,l(e_{CA})
\bigr).
$$
Together with the two displayed orders, this gives the auxiliary IET
$T_{\mathrm{aux}}$ explicitly.
\end{example}

If the refined presentation of $T_{\mathrm{aux}}$ has removable discontinuities,
we merge the corresponding adjacent intervals and denote the resulting reduced
IET by
$$
\widehat T_{\mathrm{aux}}.
$$
If there are no removable discontinuities, then
$$
\widehat T_{\mathrm{aux}}=T_{\mathrm{aux}}.
$$

\begin{definition}
\label{def:auxiliary_admissible}
We say that the ordered graph data
$$
(G,l,\{\sigma_v\}_{v\in V},\pi_0)
$$
are \emph{auxiliary-admissible} if the reduced auxiliary IET
$\widehat T_{\mathrm{aux}}$ is irreducible, satisfies Keane's i.d.o.c.~\cite{Keane}, and is
uniquely ergodic. Here ``reduced'' means that removable discontinuities in the
refined IET presentation are merged before applying irreducibility and Keane's
condition.
\end{definition}

\begin{definition}[IET-to-birth-time transfer property]
\label{def:iet_birth_transfer}
Let $T_{\mathrm{aux}}$ be the auxiliary IET associated with ordered graph
data. We say that $T_{\mathrm{aux}}$ satisfies an
$\varepsilon$-birth-time transfer estimate with function $H_T(\varepsilon)$
if, for every sufficiently small $\varepsilon>0$, every oriented edge
$a=(u\to w)$, every
$$
t\ge H_T(\varepsilon)+l(a),
$$
and every
$$
s\in[0,l(a)-\varepsilon],
$$
there exists a birth time $\tau\in\mathcal A_u$ such that
$$
0\le t-s-\tau\le\varepsilon.
$$
\end{definition}

\begin{proposition}[Conditional transfer from an auxiliary IET to saturation]
\label{prop:conditional_iet_transfer}
Assume that the auxiliary IET $T_{\mathrm{aux}}$ satisfies the
$\varepsilon$-birth-time transfer estimate of
Definition~\ref{def:iet_birth_transfer}. Then, for all sufficiently small
$\varepsilon>0$,
$$
\tau_s(\varepsilon)
\le
H_T(\varepsilon)+\max_{e\in E}l(e).
$$
\end{proposition}

\begin{proof}
Let
$$
\ell_{\max}=\max_{e\in E}l(e)
$$
and fix
$$
t\ge H_T(\varepsilon)+\ell_{\max}.
$$
We show that the hypothesis of Lemma~\ref{lem:birth_saturation_criterion} is
satisfied.

Let $a=(u\to w)\in\overrightarrow E$. Since $t\ge H_T(\varepsilon)+l(a)$,
Definition~\ref{def:iet_birth_transfer} applies. For every
$s\in[0,l(a)-\varepsilon]$, there exists $\tau\in\mathcal A_u$ such that
$$
0\le t-s-\tau\le\varepsilon.
$$
A point born at $u$ at time $\tau$ and sent along the oriented edge $a$
has, at time $t$, coordinate
$$
t-\tau=s+(t-s-\tau)
$$
measured from $u$. Hence
$$
s\le t-\tau\le s+\varepsilon\le l(a).
$$
Thus this point is on the edge $a$ at time $t$ and lies within distance
$\varepsilon$ of the point with coordinate $s$.

For the remaining interval
$$
[l(a)-\varepsilon,l(a)],
$$
apply the preceding argument to the coordinate $s_0=l(a)-\varepsilon$. This
gives a moving point whose coordinate lies in
$$
[l(a)-\varepsilon,l(a)],
$$
and hence this point is within distance $\varepsilon$ of every point in the
remaining interval.

Therefore $X_a(t)$ is an $\varepsilon$-net in $[0,l(a)]$ for every
oriented edge $a$. By Lemma~\ref{lem:birth_saturation_criterion}, the graph is
$\varepsilon$-saturated at time $t$. Since $t\ge
H_T(\varepsilon)+\ell_{\max}$ was arbitrary, the claimed bound for
$\tau_s(\varepsilon)$ follows.
\end{proof}

\subsection{Ergodicity and minimality of the auxiliary IET}
\label{subsec:ergodicity_minimality_aux}

We now record the standard qualitative dynamical properties of the auxiliary
interval exchange transformation constructed above. These properties concern
the auxiliary IET itself. By themselves they do not imply permanent saturation
of the branched graph dynamics; for saturation one also needs a birth-time
transfer statement of the form given in
Definition~\ref{def:iet_birth_transfer}.

Let
$$
(G,l,\{\sigma_v\}_{v\in V},\pi_0)
$$
be ordered graph data, and let
$$
T_{\mathrm{aux}}=T_{\pi,\lambda}
$$
be the associated auxiliary interval exchange transformation on
$$
I=\left[0,\sum_{a\in\overrightarrow E}\lambda_a\right),
\qquad
\lambda_a=l(a).
$$
The alphabet of $T_{\mathrm{aux}}$ is $\overrightarrow E$, so the number of
exchanged intervals before reduction is
$$
d=|\overrightarrow E|=2|E|.
$$

\begin{theorem}
\label{ergod}
Let $\widehat T_{\mathrm{aux}}$ be the reduced auxiliary IET associated with
the ordered graph data
$$
(G,l,\{\sigma_v\}_{v\in V},\pi_0).
$$
Assume that the permutation data of $\widehat T_{\mathrm{aux}}$ are
irreducible and that $\widehat T_{\mathrm{aux}}$ satisfies Keane's i.d.o.c.
Then $\widehat T_{\mathrm{aux}}$ is minimal.

If, in addition, the ordered graph data are auxiliary-admissible in the sense of
Definition~\ref{def:auxiliary_admissible}, then
$\widehat T_{\mathrm{aux}}$ is uniquely ergodic. Consequently, for every
continuous function $f\colon I\to\mathbb R$ and every point $x\in I$ whose
forward orbit avoids the discontinuity set of $\widehat T_{\mathrm{aux}}$, one
has
$$
\frac1N\sum_{n=0}^{N-1} f(\widehat T_{\mathrm{aux}}^n x)
\longrightarrow
\frac{1}{|I|}\int_I f(u)\,du
\qquad\text{as }N\to\infty.
$$
\end{theorem}

\begin{proof}
After removing removable discontinuities, $\widehat T_{\mathrm{aux}}$ is a
classical interval exchange transformation. By Keane's minimality criterion,
irreducibility together with the infinite distinct orbit condition implies
minimality.

If the ordered graph data are auxiliary-admissible, then, by
Definition~\ref{def:auxiliary_admissible}, the reduced auxiliary IET is uniquely
ergodic. The stated convergence is the standard equidistribution consequence of
unique ergodicity for interval exchange transformations, with respect to
normalized Lebesgue measure on $I$.
\end{proof}

\begin{remark}
Theorem~\ref{ergod} is a statement about the auxiliary IET only. In particular,
minimality or unique ergodicity of $\widehat T_{\mathrm{aux}}$ does not, by
itself, imply permanent $\varepsilon$-saturation of the branched graph
dynamics. The passage from the auxiliary IET to saturation must go through
birth-time density. In this paper this passage is made rigorously in
Theorem~\ref{main} for the rotation IET arising from two closed walks of
irrational length ratio. For a general auxiliary IET, the necessary additional
assumption is the $\varepsilon$-birth-time transfer estimate of
Definition~\ref{def:iet_birth_transfer}; under that assumption,
Corollary~\ref{cor:conditional_transfer_power} gives the corresponding
saturation bound.
\end{remark}

\begin{remark}[Rauzy--Veech estimates for the auxiliary IET]
One may apply Rauzy--Veech induction to the reduced auxiliary IET
$\widehat T_{\mathrm{aux}}$, provided its permutation data are irreducible.
If $Z(n)$ denotes the corresponding Rauzy--Veech matrix and
$$
Q(n)=Z(1)\cdots Z(n),
$$
then the growth of $Q(n)$ controls standard quantitative recurrence and
renormalization properties of $\widehat T_{\mathrm{aux}}$. Such estimates,
however, become estimates for the saturation time of the original graph
dynamics only after they have been converted into an
$\varepsilon$-birth-time transfer estimate $H_T(\varepsilon)$ in the sense
of Definition~\ref{def:iet_birth_transfer}. Establishing this conversion for
general ordered graph data remains a separate problem.
\end{remark}

\subsection{Beyond rotations: non-rotation auxiliary IETs}
\label{subsec:beyond_rotations}

The rotation construction above gives a rigorous and explicit mechanism for
permanent saturation. However, it is not the only possible source of
self-similar graph-induced interval exchange transformations. For more general
choices of the cyclic orders and of the top order $\pi_0$, the associated
auxiliary IET need not reduce to a rotation.

In such cases, one may try to find periodic Rauzy--Veech paths inside the
graph-constrained family of length vectors
$$
\lambda_a=\lambda_{\bar a}.
$$
A periodic Rauzy path whose matrix preserves these linear restrictions gives a
self-similar representative: the Perron--Frobenius eigenvector of the induced
matrix determines the corresponding edge lengths. Under the usual
irreducibility and admissibility conditions for the periodic Rauzy path,
primitivity of the corresponding matrix gives a minimal and uniquely ergodic
self-similar IET.

This point of view is also related to the study of interval exchange
transformations with linear restrictions on the length parameters, in particular
to the stability results of Dynnikov and Skripchenko for minimal and uniquely
ergodic IETs under poor restrictions~\cite{DynSkr}.

This mechanism does occur, even outside the rotation-type case. For example,
on the path graph with two edges one can choose top and bottom orders which are
not cyclic shifts of one another. A periodic Rauzy path then gives a primitive
matrix on the two edge lengths, with positive Perron--Frobenius eigenvector,
and the corresponding SAF invariant is non-zero. Thus one obtains a
graph-induced non-rotation self-similar IET with non-zero SAF.

This example is very small and should be regarded as a small model. Finding
larger non-rotation graph-induced self-similar examples with non-zero SAF,
especially on simple graphs with richer topology, appears to be a separate
interesting problem. The examples found so far suggest that this is a natural
direction for further study, but detecting such cases requires searching
simultaneously over graph combinatorics, cyclic orders, Rauzy paths, and
restricted length vectors.

\subsection{Conditional Rauzy--Veech refinement}
\label{subsec:conditional_rv_refinement}

We finally record how quantitative recurrence estimates for the auxiliary IET
would imply saturation estimates, once they have been converted into a
birth-time transfer estimate. The Rauzy--Veech mechanism is therefore used here
only conditionally: the required input is a bound on the transfer function
$H_T(\varepsilon)$ from Definition~\ref{def:iet_birth_transfer}.

\begin{corollary}[Conditional saturation estimate from a transfer bound]
\label{cor:conditional_transfer_power}
Let $G=(V,E,l)$ be a finite connected metric graph, let $v^*\in V$, and fix
ordered graph data
$$
(G,l,\{\sigma_v\}_{v\in V},\pi_0).
$$
Let $T_{\mathrm{aux}}$ be the associated auxiliary IET. Suppose that
$T_{\mathrm{aux}}$ satisfies the $\varepsilon$-birth-time transfer estimate of
Definition~\ref{def:iet_birth_transfer} with a function $H_T(\varepsilon)$.
If, for some $\kappa>0$ and $C>0$,
$$
H_T(\varepsilon)\le C\varepsilon^{-\kappa}
$$
for all sufficiently small $\varepsilon>0$, then there exists $K>0$ such that
$$
\tau_s(\varepsilon)\le K\varepsilon^{-\kappa}
$$
for all sufficiently small $\varepsilon>0$.
\end{corollary}

\begin{proof}
By Proposition~\ref{prop:conditional_iet_transfer},
$$
\tau_s(\varepsilon)
\le
H_T(\varepsilon)+\max_{e\in E}l(e).
$$
Using the assumed bound on $H_T(\varepsilon)$ gives
$$
\tau_s(\varepsilon)
\le
C\varepsilon^{-\kappa}+\max_{e\in E}l(e).
$$
The second term is independent of $\varepsilon$ and can be absorbed into the
constant for sufficiently small $\varepsilon$.
\end{proof}

\begin{corollary}[Conditional Rauzy--Veech power bound]
\label{cor:conditional_rv_power_bound}
Suppose that the hypotheses of
Corollary~\ref{cor:conditional_transfer_power} hold and that the transfer
function satisfies
$$
H_T(\varepsilon)
\le
C_\eta \varepsilon^{-\frac{1+\eta}{1-\eta}}
$$
for every $\eta>0$ and all sufficiently small $\varepsilon>0$. Then
$$
\tau_s(\varepsilon)
\le
K_\eta \varepsilon^{-\frac{1+\eta}{1-\eta}}
$$
for all sufficiently small $\varepsilon>0$.
\end{corollary}

\begin{proof}
This is Corollary~\ref{cor:conditional_transfer_power} with
$$
\kappa=\frac{1+\eta}{1-\eta}.
$$
\end{proof}

\section{Self-similar graph-induced IETs and the SAF invariant}
\label{sec:self_similar_examples}

In this section we collect several examples illustrating the auxiliary-IET
construction from Section~\ref{sec:graph_to_iet}. These examples do not enter
the proof of the rotation-based saturation theorem. Their purpose is instead to
show that the restricted families of interval exchange transformations arising
from metric graphs contain explicit self-similar representatives and exhibit
both rotation-type and non-rotation behaviour.

For every unoriented edge of a graph, the two corresponding oriented intervals
have the same length. Thus the graph construction imposes the linear
restrictions
$$
\lambda_a=\lambda_{\bar a},
$$
where $a$ and $\bar a$ are the two orientations of the same edge. This places
the auxiliary IETs in restricted affine subspaces of the full length simplex.
This point of view is related to the study of interval exchange transformations
with restrictions, in particular to the stability results of Dynnikov and
Skripchenko for minimal and uniquely ergodic IETs under poor restrictions
\cite{DynSkr}.

For an interval exchange transformation $T$ with interval lengths
$\lambda_\alpha$ and translation lengths $\delta_\alpha$, the
Sah--Arnoux--Fathi invariant is
$$
\operatorname{SAF}(T)
=
\sum_\alpha \lambda_\alpha\wedge_{\mathbb Q}\delta_\alpha
\in
\mathbb R\wedge_{\mathbb Q}\mathbb R.
$$
Here $\delta_\alpha$ is defined by
$$
T(x)=x+\delta_\alpha,\qquad x\in I_\alpha.
$$
The exterior product is taken over $\mathbb Q$, so it is bilinear over
$\mathbb Q$ and satisfies
$$
x\wedge_{\mathbb Q}y=-y\wedge_{\mathbb Q}x,
\qquad
x\wedge_{\mathbb Q}x=0.
$$

The simplest examples are rotation-type. In these examples the cyclic orders
are chosen so that the successor map
$$
\Phi(a)=\sigma_{h(a)}(\bar a)
$$
is one cycle on the set of oriented edges. If the top order is chosen along this
cycle, then the bottom order is a cyclic shift, and the associated auxiliary IET
reduces to a circle rotation. Choosing a quadratic irrational rotation number
then gives a self-similar minimal and uniquely ergodic IET. Since the rotation
number is irrational, the Sah--Arnoux--Fathi invariant is non-zero.

This construction gives examples on several simple graphs, including two
triangles sharing a vertex, the complete bipartite graphs $K_{2,3}$ and
$K_{3,3}$, and the complete graph $K_5$. These examples are dynamically simple,
because the resulting IETs are rotations, but they show that graph-induced
restricted families naturally contain explicit self-similar IETs with non-zero
SAF invariant.

The first non-rotation example is a small model example on the path graph with two
edges,
$$
v_0 \stackrel{e_1}{\longleftrightarrow} v_1
\stackrel{e_2}{\longleftrightarrow} v_2.
$$
Let
$$
a=e_1^+=(v_0\to v_1),\qquad \bar a=e_1^-=(v_1\to v_0),
$$
and
$$
b=e_2^+=(v_1\to v_2),\qquad \bar b=e_2^-=(v_2\to v_1).
$$
Choose cyclic orders
$$
\sigma_{v_0}=(a),\qquad
\sigma_{v_1}=(\bar a,\ b),\qquad
\sigma_{v_2}=(\bar b).
$$
Then
$$
\Phi(a)=b,\qquad
\Phi(\bar a)=a,\qquad
\Phi(b)=\bar b,\qquad
\Phi(\bar b)=\bar a.
$$
For the top order
$$
\pi_0=(a,\bar a,b,\bar b)
$$
the corresponding bottom order is
$$
\pi_1=(b,a,\bar b,\bar a).
$$
These two orders are not cyclic shifts of one another, so the auxiliary IET is
not rotation-type.

A periodic Rauzy path for this IET is
$$
\gamma=01001101010.
$$
The induced matrix on the two unoriented edge lengths is
$$
A=
\begin{pmatrix}
4&3\\
5&4
\end{pmatrix}.
$$
This matrix is primitive. Its Perron--Frobenius eigenvalue is
$$
\rho=4+\sqrt{15},
$$
and a positive eigenvector is
$$
\ell=
\begin{pmatrix}
1\\
\sqrt{15}/3
\end{pmatrix}.
$$
Indeed,
$$
A
\begin{pmatrix}
1\\
\sqrt{15}/3
\end{pmatrix}
=
(4+\sqrt{15})
\begin{pmatrix}
1\\
\sqrt{15}/3
\end{pmatrix}.
$$
Therefore the choice
$$
l(e_1)=1,\qquad l(e_2)=\frac{\sqrt{15}}{3}
$$
gives a self-similar graph-induced non-rotation IET.

For the generic paired length family
$$
\lambda_a=\lambda_{\bar a}=x,\qquad
\lambda_b=\lambda_{\bar b}=y,
$$
one computes
$$
\operatorname{SAF}(T)=6\,x\wedge_{\mathbb Q}y.
$$
For the Perron--Frobenius lengths above this becomes
$$
\operatorname{SAF}(T)
=
6\,1\wedge_{\mathbb Q}\frac{\sqrt{15}}{3}
=
2\,1\wedge_{\mathbb Q}\sqrt{15}
\neq0.
$$
Thus this example is graph-induced, self-similar, non-rotation, and has
non-zero SAF invariant. The underlying graph is very small, so the example
should be viewed as a small model rather than as a structurally rich graph
example.

The ordered $K_4$ example used in the numerical section gives a different
phenomenon. For the edge lengths
$$
\sqrt2,\sqrt3,\sqrt5,\sqrt7,\sqrt{11},\sqrt{13},
$$
and for the cyclic orders specified in Section~\ref{sec:numerics}, a direct
computation gives a non-zero SAF invariant. Thus this restricted $K_4$ family
is a natural poor-type candidate in the sense of the theory of IETs with
restrictions. However, we have not found a self-similar representative inside
this $K_4$ family.

Finally, the Arnoux--Yoccoz/Arnoux--Rauzy construction gives a non-rotation
self-similar example related to a theta graph or circle exchange with three
pairs of equal lengths. In the cubic case, if
$$
\alpha^3+\alpha^2+\alpha=1,
$$
the corresponding length pattern is
$$
\frac{\alpha}{2},\frac{\alpha}{2},
\frac{\alpha^2}{2},\frac{\alpha^2}{2},
\frac{\alpha^3}{2},\frac{\alpha^3}{2}.
$$
This example has known self-similarity and unique ergodicity, but its SAF
invariant is zero. Therefore it does not directly provide a poor-restriction
stability example.

The examples above suggest the following picture. Rotation-type graph-induced
examples with non-zero SAF are easy to construct once the successor permutation
is one cycle. Non-rotation self-similar examples with non-zero SAF also exist,
as shown by the two-edge path example, but finding larger examples on simple
graphs with richer topology appears to be a separate problem.

\section{Numerical illustrations and possible applications}\label{sec:numerics}

In this section we present numerical illustrations of the upper-bound expression derived above and briefly indicate possible interpretations of the model in network-type settings. The computations are intended to illustrate how the bound behaves for concrete graph families with incommensurable edge lengths.

\medskip
\noindent\textbf{Graph families and choice of parameters.}
For the numerical experiments we considered two representative graph families: the complete graph \(K_4\) and a star graph with five leaves. These examples were chosen because they exhibit different branching and connectivity patterns while remaining simple enough to simulate explicitly. In each case the edge lengths were chosen to be incommensurable:
\begin{itemize}
    \item for the graph \(K_4\), the edge lengths were taken to be
    \[
    \sqrt{2},\sqrt{3},\sqrt{5},\sqrt{7},\sqrt{11},\sqrt{13};
    \]
    \item for the star graph with five leaves, the edge lengths were taken to be
    \[
    \sqrt{2},\sqrt{3},\sqrt{5},\sqrt{7},\sqrt{11}.
    \]
\end{itemize}

\medskip
\noindent\textbf{Cyclic orders used in the simulations.}
For the \(K_4\) example, write the vertices as \(v_1,v_2,v_3,v_4\), and denote the oriented edge from \(v_i\) to \(v_j\) by
\[
a_{ij}=(v_i\to v_j).
\]
We used the cyclic orders
\[
\sigma_{v_1}=(a_{12}\ a_{13}\ a_{14}),\qquad
\sigma_{v_2}=(a_{21}\ a_{23}\ a_{24}),
\]
\[
\sigma_{v_3}=(a_{31}\ a_{32}\ a_{34}),\qquad
\sigma_{v_4}=(a_{41}\ a_{42}\ a_{43}).
\]
The top order was chosen as
\[
\pi_0=(a_{12},a_{13},a_{14},a_{21},a_{23},a_{24},
a_{31},a_{32},a_{34},a_{41},a_{42},a_{43}).
\]
Using the rule \(\pi_1(\Phi(a))=\pi_0(a)\), this gives the bottom order
\[
\pi_1=(a_{23},a_{32},a_{42},a_{13},a_{34},a_{43},
a_{14},a_{24},a_{41},a_{12},a_{21},a_{31}).
\]
The length assignment was
\[
l(e_{12})=\sqrt2,\quad l(e_{13})=\sqrt3,\quad l(e_{14})=\sqrt5,
\]
\[
l(e_{23})=\sqrt7,\quad l(e_{24})=\sqrt{11},\quad l(e_{34})=\sqrt{13},
\]
and the oriented-edge length vector was defined by
\[
\lambda_{a_{ij}}=l(e_{ij}).
\]
For this ordered $K_4$ data, the induced permutation pair is irreducible. The
corresponding auxiliary IET is not rotation-type. Since the main theorem of the
paper is the rotation-based saturation estimate, the $K_4$ plot should be
understood as a numerical illustration of the broader auxiliary-IET framework
rather than as a direct application of Theorem~\ref{main}.

For the star graph with five leaves, write
\[
a_i=(o\to v_i),\qquad b_i=(v_i\to o),\qquad i=1,\dots,5.
\]
We used the cyclic order
\[
a_1\mapsto a_2\mapsto a_3\mapsto a_4\mapsto a_5\mapsto a_1
\]
at the center and the trivial order at each leaf. With the top order
\[
\pi_0=(a_1,b_1,a_2,b_2,a_3,b_3,a_4,b_4,a_5,b_5),
\]
the associated IET is a rotation on an interval of length
\[
L=2(\ell_1+\ell_2+\ell_3+\ell_4+\ell_5).
\]
It is minimal and uniquely ergodic whenever
\[
\frac{\ell_1}{L}\notin\mathbb Q.
\]
For the square-root lengths used above, this condition holds. Thus the star
graph gives a direct illustration of the rotation mechanism used in
Theorem~\ref{main}.

\medskip
\noindent\textbf{Simulation procedure.}
For each graph and each value of $\varepsilon$, we simulated the dispersing
moving-point dynamics and recorded the moment of permanent
$\varepsilon$-saturation. For the star graph, we compared these values with the
rotation-based upper-bound expression from Theorem~\ref{main}. For $K_4$, the
comparison should be interpreted as a numerical illustration of the auxiliary
IET perspective for the ordered graph data specified above.

\medskip
\noindent\textbf{Observed behaviour.}
For both graph families, the observed moments of permanent
$\varepsilon$-saturation remain below the plotted upper-bound expression
throughout the tested range of $\varepsilon$. For the star graph, this gives a
direct numerical illustration of the rotation mechanism used in
Theorem~\ref{main}. For $K_4$, the computation should be understood as a
numerical illustration of the broader auxiliary-IET framework for the ordered
graph data specified above, rather than as a direct application of
Theorem~\ref{main}. In both cases, the visible gap between the observed values
and the bound suggests that the estimate is conservative and may admit
refinement under additional geometric or dynamical assumptions.

\begin{figure}[t]
    \centering
    \subfloat[\centering Graph \(K_4\)]{{\includegraphics[scale=0.26]{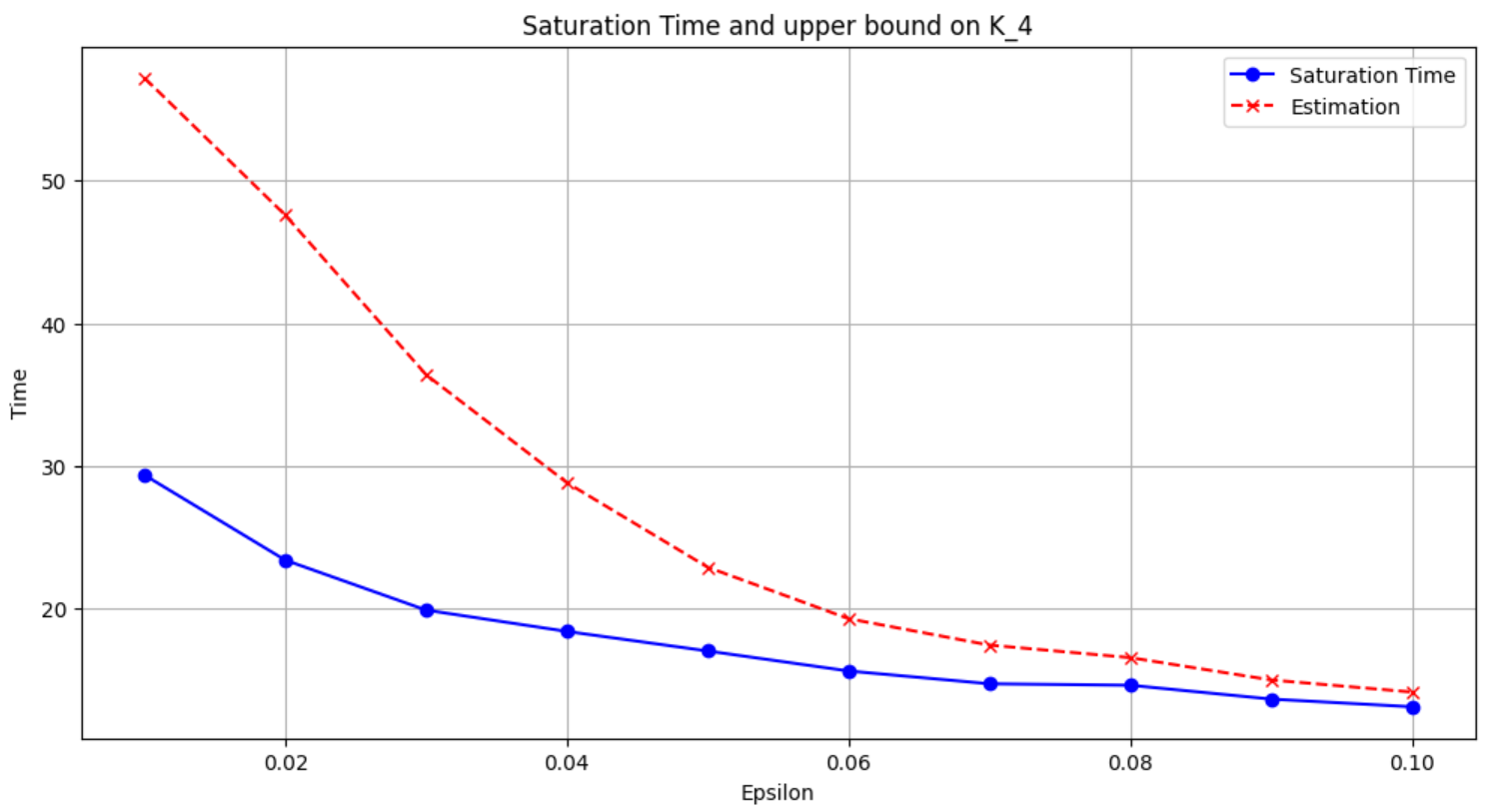} }}%
    \hfill
    \subfloat[\centering Star graph]{{\includegraphics[scale=0.26]{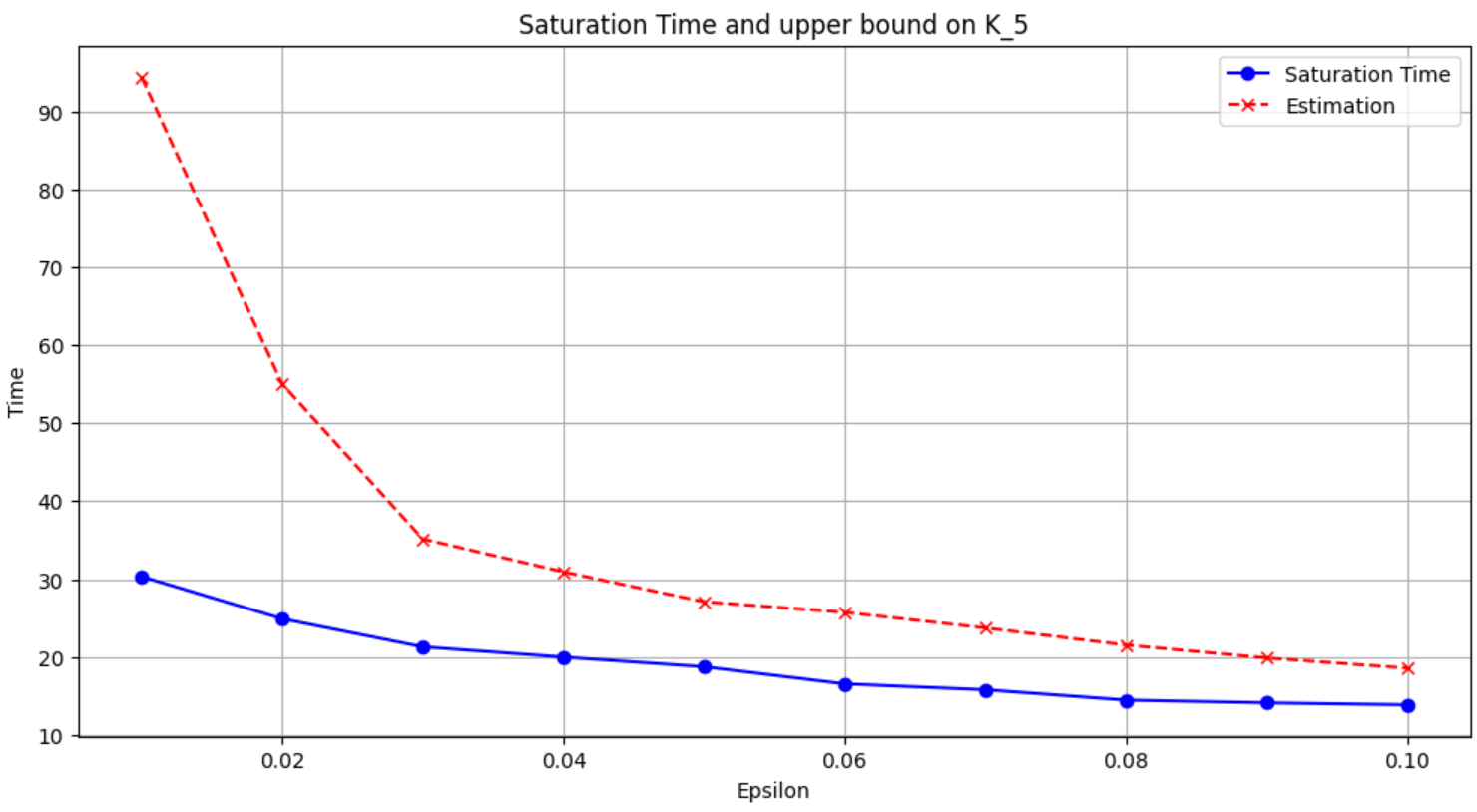} }}%
    \caption{Observed moments of permanent $\varepsilon$-saturation (solid blue) and the plotted upper-bound expression (dashed red).}
    \label{fig:example}
\end{figure}

\medskip
\noindent\textbf{Possible application interpretation.}
Although our results are purely mathematical, the model may be viewed as an idealized description of branching transport on a finite network. For example, one may think of a signal, particle, or packet that travels along edges and, upon reaching a vertex, generates copies that continue along other admissible edges. In that interpretation, the moment of permanent \(\varepsilon\)-saturation gives a quantitative bound on the time after which the network is persistently covered at spatial resolution \(\varepsilon\). From this perspective, the theorem provides a conservative estimate for the onset of persistent coverage in a finite branching network with incommensurable travel times.

\medskip
\noindent\textbf{Limitations.}
The present model is deterministic and highly idealized: it does not include capacities, stochastic delays, dissipation, queueing effects, or optimization constraints. Consequently, the discussion above should be understood only as a possible interpretation of the mathematical framework, rather than as a direct application to an engineering system.

\section{Discussion and Open Questions}\label{sec:discussion}

In this work we studied the moment of permanent $\varepsilon$-saturation for
branched moving-point dynamics on finite metric graphs. The main rigorous
result is the rotation-based estimate of Theorem~\ref{main}. Its proof is based
on the birth-time formulation of the dynamics and on the following implication:
if two closed walks based at the initial vertex have incommensurable lengths
$A$ and $B$, then the arithmetic semigroup
$$
\Gamma_{A,B}=\{nA+mB:n,m\in\mathbb Z_{\ge0}\}
$$
is sufficiently dense at large times to force permanent saturation.

The auxiliary oriented-edge IET constructed in Section~\ref{sec:graph_to_iet}
plays a different role. It encodes the ordered edge-state combinatorics of the
graph and provides a natural setting for Rauzy--Veech renormalization. However,
it is not a complete model of the branched dynamics, since the latter is
non-invertible, whereas an IET is invertible away from its discontinuities.
Thus the general mechanism suggested by the auxiliary IET has the form
\[
\begin{aligned}
\text{recurrence/discrepancy of } \widehat T_{\mathrm{aux}}
&\;\Longrightarrow\; \text{birth-time density} \\
&\;\Longrightarrow\; \text{permanent saturation}.
\end{aligned}
\]
The second implication is given by Proposition~\ref{prop:conditional_iet_transfer};
the first one is precisely the additional transfer problem.

\subsection{Topological interpretation}

Let $T_{\mathrm{aux}}=T_{\pi,\lambda}$ be the auxiliary IET associated with
ordered graph data
$$
(G,l,\{\sigma_v\}_{v\in V},\pi_0).
$$
As a classical IET, $T_{\mathrm{aux}}$ admits a suspension to a translation
surface. If the reduced IET has $d$ exchanged intervals and permutation data
$\pi$, the corresponding suspension carries the usual skew-symmetric
intersection form associated with Rauzy--Veech induction. In the classical
case, the rank of this form determines the genus $g$ of the suspension surface:
$$
\operatorname{rank}\Omega_\pi=2g.
$$
This is a standard part of the relation between IETs, suspensions, and
translation surfaces; see \cite{Ve2,ViaT,AvGoYo}.

For graph-induced auxiliary IETs, the same classical construction applies to
$\widehat T_{\mathrm{aux}}$. What is not yet clear is whether the resulting
translation surface has an intrinsic interpretation in terms of the original
branched graph dynamics. The obstruction is structural: the suspension flow of
an IET is invertible away from singularities, whereas the graph dynamics
contains branching and multiplicities. Thus the suspension surface is naturally
attached to the auxiliary IET, but not automatically to the full branched
system.

This leads to the following problem.

\begin{problem}
Given ordered graph data
$$
(G,l,\{\sigma_v\}_{v\in V},\pi_0),
$$
describe the topology of the suspension surface of the reduced auxiliary IET
$\widehat T_{\mathrm{aux}}$ directly in terms of the graph $G$ and the cyclic
orders $\{\sigma_v\}_{v\in V}$. In particular, determine whether the genus
$$
g=\frac12\operatorname{rank}\Omega_\pi
$$
admits an intrinsic graph-theoretic interpretation.
\end{problem}

\subsection{Symbolic complexity and renormalization complexity}

There are two different notions of complexity associated with an IET. The first
one is symbolic complexity of the natural coding. Let $T$ be an IET on $d$
subintervals, and let $p_T(n)$ denote the number of words of length $n$ in the
natural coding by the top partition:
$$
p_T(n)=
\#\left\{
\alpha_0\alpha_1\cdots\alpha_{n-1}:
\exists x\in I,\ T^k x\in I_{\alpha_k},\ k=0,\ldots,n-1
\right\}.
$$

\begin{proposition}
Let $T$ be a minimal IET on $d$ intervals satisfying Keane's i.d.o.c. Then
$$
p_T(n)=(d-1)n+1
\qquad
\text{for all } n\ge1.
$$
\end{proposition}

\begin{proof}
For a regular IET satisfying Keane's i.d.o.c., the natural coding has exactly
$d-1$ right-special words of each length. Hence
$$
p_T(n+1)-p_T(n)=d-1
\qquad
\text{for all } n\ge1.
$$
Since $p_T(1)=d$, summing the differences gives
$$
p_T(n)=d+(n-1)(d-1)=(d-1)n+1.
$$
This is the classical complexity formula for regular interval exchange
transformations; see \cite{Keane,Vian}.
\end{proof}

For the auxiliary IET associated with a metric graph before reduction, the
alphabet is $\overrightarrow E$, and hence
$$
d=|\overrightarrow E|=2|E|.
$$
If the reduced auxiliary IET is minimal and satisfies Keane's i.d.o.c., its
natural coding therefore has linear symbolic complexity. However, this fact
does not by itself imply any quantitative saturation estimate.

The second notion is renormalization complexity. It is encoded by the growth of
Rauzy--Veech matrices. If $Z(n)$ denotes the matrix of the $n$th Rauzy--Veech
step and
$$
Q(n)=Z(1)Z(2)\cdots Z(n),
$$
then quantitative recurrence and discrepancy estimates are governed by the
growth of $\|Q(n)\|$, by balance properties of the corresponding towers, and,
in the suspension picture, by Lyapunov exponents of the Kontsevich--Zorich
cocycle. Thus symbolic complexity gives the combinatorial size of the language,
whereas renormalization complexity controls the quantitative distribution of
orbits.

In the graph setting, a saturation estimate must involve the latter type of
information. More precisely, one needs a chain of estimates of the form
\[
\begin{aligned}
\|Q(n)\|\ \text{controlled}
&\;\Longrightarrow\; D_N(x;\widehat T_{\mathrm{aux}})\ \text{controlled} \\
&\;\Longrightarrow\; H_T(\varepsilon)\ \text{controlled} \\
&\;\Longrightarrow\; \tau_s(\varepsilon)\ \text{controlled}.
\end{aligned}
\]
The final implication is Proposition~\ref{prop:conditional_iet_transfer}. The
middle implication, from discrepancy of the auxiliary IET to birth-time
density, is the main missing step in the general non-rotation case.

\subsection{Conjectural Lyapunov-type refinements}

Let $T=T_{\pi,\lambda}:I\to I$ be the reduced auxiliary IET associated with
ordered graph data, and let $L=|I|$. For $x\in I$ and $N\ge1$, define the
interval discrepancy of the orbit segment
$$
x,\ Tx,\ \ldots,\ T^{N-1}x
$$
by
$$
D_N(x;T)
=
\sup_{J\subset I}
\left|
\frac{1}{N}\sum_{k=0}^{N-1}\mathbf 1_J(T^k x)
-
\frac{|J|}{L}
\right|,
$$
where the supremum is taken over all intervals $J\subset I$.

Classical discrepancy estimates for IETs and translation flows suggest bounds
of the form
$$
D_N(x;T)
\le
C_T\frac{(\log N)^m}{N},
$$
where the constant $C_T$ is governed by the Rauzy--Veech or
Kontsevich--Zorich renormalization data. If such a discrepancy estimate can be
converted into an $\varepsilon$-birth-time transfer estimate, then
Corollary~\ref{cor:conditional_transfer_power} gives a corresponding saturation
bound.

This motivates the following conjecture.

\begin{conjecture}
\label{question_conj}
Let $G=(V,E,l)$ be a connected metric graph with chosen cyclic orders, and let
$T=T_{\pi,\lambda}$ be the corresponding reduced auxiliary IET. Suppose that the
ordered graph data are auxiliary-admissible and that the relevant
renormalization cocycle has non-negative Lyapunov exponents
$$
\chi_1\ge \chi_2\ge\cdots\ge \chi_m\ge0.
$$
Assume moreover that discrepancy estimates for $T$ imply the
$\varepsilon$-birth-time transfer property of
Definition~\ref{def:iet_birth_transfer}. Then one expects an upper bound of the
form
$$
\tau_s(\varepsilon)
\le
\frac{K\exp\left(\sum_{i=1}^{m}\chi_i\right)}{\varepsilon}
\left(
\log\left(K\exp\left(\sum_{i=1}^{m}\chi_i\right)\right)
-\log\varepsilon
\right)^m,
$$
where $K>0$ depends on the graph data.
\end{conjecture}

The logarithmic factor is consistent with classical discrepancy estimates for
IETs and translation flows, where logarithmic losses naturally appear in
quantitative uniform distribution estimates~\cite{AvGoYo,Disrep,Saussl}. The
conjecture is not proved in the present paper. Its main unproved component is
not the classical IET discrepancy estimate, but the transfer from discrepancy
of $\widehat T_{\mathrm{aux}}$ to density of the birth-time sets
$\mathcal A_v$.

Equivalently, the desired missing implication can be formulated as follows:
find explicit constants $C>0$ and $\kappa>0$ such that
$$
D_N(x;\widehat T_{\mathrm{aux}})
\le
C N^{-1}(\log N)^m
$$
implies
$$
H_T(\varepsilon)\le C'\varepsilon^{-\kappa}
|\log\varepsilon|^{m'}
$$
for some constants $C'>0$ and $m'\ge0$ depending on the graph data. Then
Proposition~\ref{prop:conditional_iet_transfer} would give
$$
\tau_s(\varepsilon)
\le
C'\varepsilon^{-\kappa}
|\log\varepsilon|^{m'}
+
\max_{e\in E}l(e).
$$

\begin{figure}%
    \centering
    \subfloat[\centering Graph $K_4$]{{\includegraphics[scale=0.26]{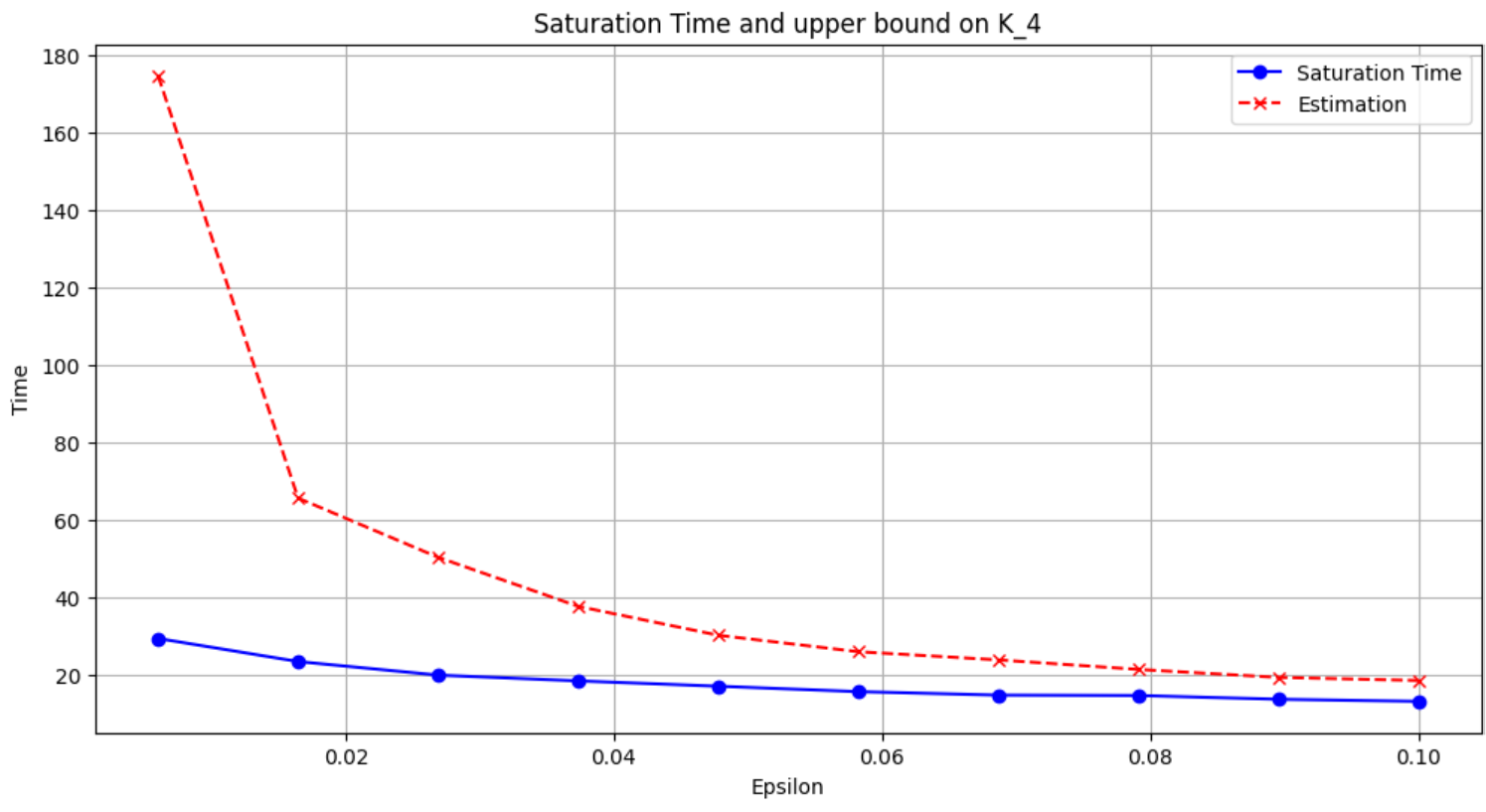} }}%
    \hfill
    \subfloat[\centering Star graph]{{\includegraphics[scale=0.26]{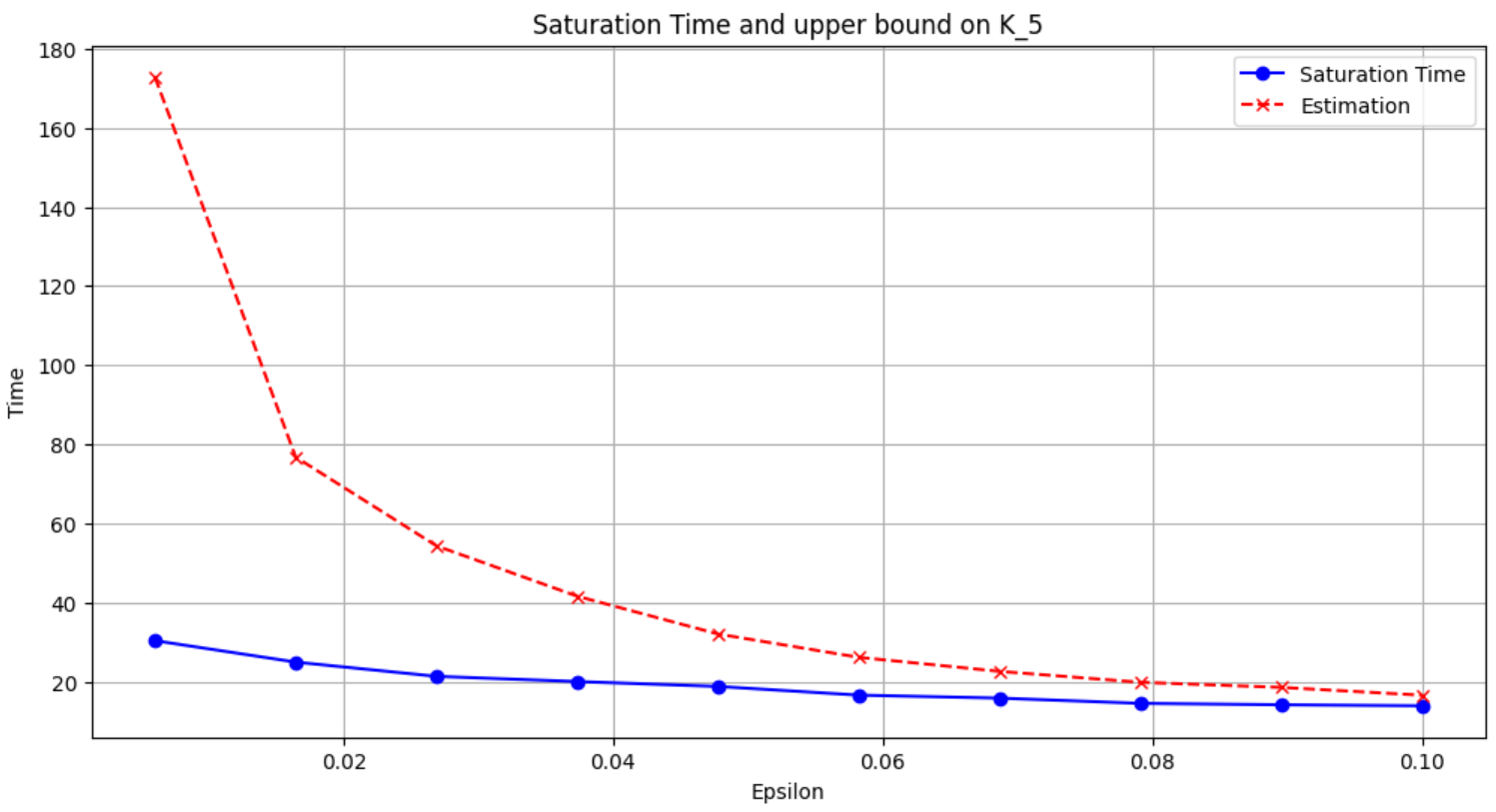} }}%
    \caption{Numerical comparison with the conjectural Lyapunov-type estimate.}%
    \label{fig:props}%
\end{figure}

We compared the conjectural expression in
Conjecture~\ref{question_conj} with numerical data for the graph families
considered in Section~\ref{sec:numerics}; see Figure~\ref{fig:props}. In the
displayed examples we set $K=1$.

The estimated curves remain above the observed moments of permanent
$\varepsilon$-saturation in the examples considered. These computations are
illustrative only and should not be interpreted as evidence for the missing
IET-to-birth-time transfer theorem.

\subsection{Open problems}

We end by recording several open problems suggested by the preceding discussion.

\begin{problem}[IET-to-birth-time transfer]
Find sufficient conditions on ordered graph data
$$
(G,l,\{\sigma_v\}_{v\in V},\pi_0)
$$
under which recurrence or discrepancy estimates for the reduced auxiliary IET
$\widehat T_{\mathrm{aux}}$ imply an $\varepsilon$-birth-time transfer estimate
in the sense of Definition~\ref{def:iet_birth_transfer}. In particular,
determine whether bounds of the form
$$
D_N(x;\widehat T_{\mathrm{aux}})
\le
C N^{-1}(\log N)^m
$$
imply explicit estimates for $H_T(\varepsilon)$.
\end{problem}

\begin{problem}[Auxiliary admissibility]
Characterize the ordered graph data for which the reduced auxiliary IET
$\widehat T_{\mathrm{aux}}$ is irreducible, satisfies Keane's i.d.o.c., and is
uniquely ergodic. Since graph-induced length vectors satisfy
$$
\lambda_a=\lambda_{\bar a},
$$
this problem naturally belongs to the theory of interval exchange
transformations with restrictions, in the sense of Dynnikov and
Skripchenko~\cite{DynSkr}.
\end{problem}

\begin{problem}[Non-rotation self-similar examples]
Find larger simple graphs whose graph-induced restricted families contain
non-rotation self-similar IETs with non-zero Sah--Arnoux--Fathi invariant. The
small path example in Section~\ref{sec:self_similar_examples} shows that such
phenomena occur, while the rotation-type examples show that non-zero SAF is
easy to obtain when the successor permutation is one cycle. The problem is to
find structurally richer examples satisfying
$$
\text{non-rotation}
\quad+\quad
\text{self-similarity}
\quad+\quad
\operatorname{SAF}(T)\neq0.
$$
\end{problem}

\begin{problem}[Topological realization]
Determine whether the suspension surface of the reduced auxiliary IET has an
intrinsic interpretation in terms of the original metric graph dynamics. In
particular, describe the meaning of
$$
g=\frac12\operatorname{rank}\Omega_\pi
$$
for graph-induced auxiliary IETs and determine whether it can be computed
directly from the graph and the cyclic orders.
\end{problem}

\begin{problem}[Extensions of the model]
Extend the theory to variants of the branched dynamics, such as directed graphs,
weighted branching rules, non-unit velocities, stochastic delays, or metric
graphs with capacities and dissipation. In each case, one must replace the
birth-time arithmetic used in Theorem~\ref{main} by an appropriate analogue.
\end{problem}

\section*{Declarations}

\subsection*{Funding}
This research received no specific grant from any funding agency, commercial or not-for-profit sectors.

\section*{Acknowledgments}

We thank Sergei Nechaev for helpful discussions. The article was prepared within the framework of the HSE University Basic Research Program as part of regular institutional support; no specific grant was received.

\bibliographystyle{amsplain}

\end{document}